\documentstyle{article}
\title{Implicit Function Theorem for systems of polynomial equations
with vanishing Jacobian and its application to flexible polyhedra and
frameworks
\thanks{This work was partially supported by RFBR grant 98--01--00688
and by INTAS-RFBR grant IR--97--1778.}}

\author{Victor Alexandrov}

\date{Sobolev Institute of Mathematics, Novosibirsk-90, 630090, Russia.
alex@math.nsc.ru}

\begin{document}

\maketitle

\begin{abstract}
We study the existence problem for a local implicit function determined by
a system of nonlinear algebraic equations in the particular case when
the determinant of its Jacobian matrix vanishes at the point under consideration.
We present a system of sufficient conditions that implies existence of a local
implicit function as well as another system of sufficient conditions that
guarantees absence of a local implicit function.
The results obtained are applied to proving new and classical results
on flexibility and rigidity of polyhedra and frameworks.


\end{abstract}

\section{Introduction}

Let $F: \mbox{\bf R}^l\times \mbox{\bf R}^m \to \mbox{\bf R}^n$ be a
differentiable mapping, let $t,t_0\in \mbox{\bf R}^l$,
let $X,X_0\in \mbox{\bf R}^m$, and let $F(t_0,X_0)=0$.
The classical Implicit Function Theorem provides conditions which imply
that the equation $F(t,X)=0$ determines an implicit function $X=X(t)$
in a neighborhood of the point $(t_0,X_0)$.
The principal condition is invertibility of the operator $F'_X(t_0,X_0)$.

The Inverse Function Theorem has numerous applications and is generalized
in various directions.
However, the author is not aware of any version of this theorem
which guarantees existence of an implicit function in the case when the
operator $F'_X(t_0,X_0)$ is not invertible.
In the present paper we will partially fill in this gap.

Our study is motivated by that of flexible polyhedra and frameworks.
It turns out that mappings $F$ which appear is that field do not depend
on the parameter $t$.
We will focus our attention on this particular case.
The following system of nonlinear algebraic equations can be considered
as a typical example of a system to which our arguments can be applied:
$$ F_1(t,x_1,x_2,x_3)\equiv x_1^2+x_2^2-x_3^2-1 =0, $$
$$ F_2(t,x_1,x_2,x_3)\equiv 3x_1+x_2-3x_3+1 =0,\eqno(1) $$ $$
F_3(t,x_1,x_2,x_3)\equiv x_1-3x_2+x_3+3 =0. $$
The parameter $t$ is not explicitly involved in this system.
The point  $X_0=(5,5,7)^{\rm T}$ satisfies (1).
The determinant of the Jacobian matrix of (1) vanishes
at $X_0=(x_1,x,_2,x_3)$:
$$ \det F'_X(t,X_0)=
\left|
\begin{array}{ccc}
2x_1 & 2x_2 & -2x_3\\
3 & 1 & -3 \\
1 & -3 & 1
\end{array}
\right|
=
\left|
\begin{array}{ccc}
10 & 10 & -14 \\
3 & 1 & -3 \\
1 & -3 & 1
\end{array}
\right|=0.
$$
Therefore, the classical Implicit Function Theorem cannot be applied to (1).
Nevertheless, from the results presented below it will follow that $X_0$
is not an isolated solution to (1); on the contrary, it belongs to a
continuous family of solutions $X=X(t)$ which can be treated as an
implicit function determined by system (1) and initial point $X_0$.

\section{Sufficient conditions for existence of an implicit function}

Let $X=(x_1,\dots ,x_m)\in \mbox{\bf R}^m$ and let
$F(X)=(F_1(X),\dots ,F_n(X))$,
where each $F_k$ $(k=1,\dots ,n)$ is a polynomial.
Without loss of generality, we may assume that the degree
of each $F_k$ is at most 2.

To explain the last statement, we assume, for example, that
in the system $F(X)=0$ under consideration each polynomial
$F_k$ $(k=1,\dots ,n-1)$ is of degree at most 2
while $F_n(X)=x_1^2x_2-1$.
Introduce a new independent variable  $x_{m+1}$
and put $\widetilde X= (x_1,\dots ,$ $x_m, x_{m+1})$.
We also introduce new functions
$\widetilde F_n (\widetilde X)= x_{m+1}x_2-1$ and
$\widetilde F_{n+1}=x_{m+1}-x_1^2$ and put $\widetilde F
(\widetilde X)= (F_1(X), \dots , F_{n-1}(X), \widetilde
F_n(\widetilde X), \widetilde F_{n+1} (\widetilde X))$.
Obviously, the system  $F(X)=0$ is equivalent to the system
$\widetilde F(\widetilde X))=0$ and each equation of the latter
system is of degree at most 2.

Thus, without loss of generality, we may assume that the degree
of each $F_k$ is at most 2.
In this case, $F_k$ can written as
$$ F_k(X)=\sum_{i=1}^{m} \sum_{j=1}^{m} \alpha
_{ij}^kx_ix_j + \sum_{i=1}^{m}\beta_i ^k x_i + \gamma ^k, $$ where
$\alpha _{ij}^k$, $\beta_i ^k$, and  $\gamma ^k$ are some reals
satisfying $\alpha _{ij}^k=\alpha _{ji}^k$.

It is well known that, if a system of polynomial equations admits
a family of solutions which is continuous with respect to a parameter,
then this system also admits a family of solutions which depends
analytically on a parameter (possibly different); see, for example,
\cite{Gluck} or Lemma 18.3 in \cite{Wallace}.
Thus, assuming that the system  $F(X)=0$ admits a continuous family
of solutions $X=X(t)\equiv(x_1(t),\dots ,x_m(t))$,
we may assume without loss of generality that this family depends
analytically on $t$, i.e., it can be expanded in a Maclaurin
series:
$$ x_i(t)=\sum_{k=0}^{\infty} x_{i,k}
t^k, \qquad x_{i,k}\in \mbox{\bf R}. $$
Substituting this expansion into equation $F_k(X)=0$,
we obtain $$ \sum_{i=1}^{m}\sum_{j=1}^{m}
\biggl[\alpha_{ij}^{k} \biggl(\sum_{p=0}^{\infty}x_{i,p}t^p\biggr)
\biggl(\sum_{q=0}^{\infty}x_{i,q}t^q\biggr)\biggr]
+\sum_{i=1}^{m}\beta_i^k\biggl(\sum_{p=0}^{\infty}x_{i,p}t^p\biggr)
+\gamma^k =0 $$ or $$ \sum_{p=0}^{\infty} \biggl[
\sum_{i=1}^{m}\sum_{j=1}^{m}
\alpha_{ij}^{k}\sum_{q=0}^{p}x_{i,q}x_{j,p-q} \biggr]t^p +
\sum_{p=0}^{\infty} \biggl[\sum_{i=1}^{m} \beta_i^k
x_{i,p}t^p\biggr] +\gamma^k =0. $$
Interpreting the left-hand side of the latter equation as the
Maclaurin expansion of the function that equals zero identically,
we conclude that, in this expansion, the coefficient of $t^p$ equals
zero for each $p\geq 1$, i.e., the equation
$$\sum_{i=1}^{m}\sum_{j=1}^{m}\sum_{q=0}^{p}
\alpha_{ij}^{k}x_{i,q}x_{j,p-q} + \sum_{i=1}^{m} \beta_i^k x_{i,p}
+\gamma^k =0 \eqno(2) $$
holds for all $p\geq 1$ and $1\leq k\leq n$.

For each  $p\geq 1$, put $X_p=(x_{1,p},x_{2,p},\dots ,x_{m,p})\in\mbox{\bf R}^m$.
Let the bilinear mapping
$B: \mbox{\bf R}^m\times \mbox{\bf R}^m\to \mbox{\bf R}^n$
be defined by the following rule:
if $X=(x_1,\dots, x_m)\in \mbox{\bf R}^m$ and $Y=(y_1,\dots , y_m)\in \mbox{\bf R}^m$
then the $k$th coordinate of the vector $B(X,Y)$ is equal to
$$
\sum_{i=1}^{m}\sum_{j=1}^{m} \alpha_{ij}^{k}x_{i}y_{j}.
$$
Let the linear mapping
$A:\mbox{\bf R}^m\to \mbox{\bf R}^n$ be defined by the following rule:
if $X=(x_1,\dots, x_m)\in \mbox{\bf R}^m$ then
the $k$th coordinate of the vector $A(X)$ equals
$$
\sum_{i=1}^{m} \beta_i^k x_{i}.
$$
Using this notation, we can rewrite (2) as
$$
\sum_{p=0}^{q} B(X_p,X_{q-p}) + AX_q=0.
$$

It follows that, if the vectors $X_0$, $X_1$, \dots $X_{q-1}$ are given and
we seek the vector $X_q$, then we need to solve the following
system of linear equations:
$$
B(X_0,X_q)+B(X_q,X_0)+AX_q=
-\sum_{p=1}^{q-1}B(X_p,X_{q-p}). \eqno(3)
$$

Let the linear mapping $C:\mbox{\bf R}^m\to \mbox{\bf R}^n$ be given by
$CX=B(X_0,X)+B(X,X_0)+AX$.
Then we can rewrite (3) in concise form:
$$
CX_q=-\sum_{p=1}^{q-1}B(X_p,X_{q-p}). \eqno(4)
$$

In the preceding consideration, the vectors $X_p$ were generated by
the Maclaurin coefficients $x_{i,p}$ of a family of exact solutions
to $F(X)=0$.
Now we assume that we have an arbitrary finite set $Y_0$, $Y_1$, \dots , $Y_q$
of vectors in $\mbox{\bf R}^m$.
We call the expression
$$
Y(t)=\sum_{p=0}^{q} Y_pt^p
$$
an {\it approximate solution of degree $q$} to the system of polynomial
equations $F(X)=0$ if, for each $p=1,2,\dots ,q$, the coefficient
of $t^p$ in the Maclaurin expansion of the function $F(Y(t))$ is equal to zero.
An equivalent formulation of this condition is as follows:
for each $p=1,2,\dots ,q$, the equation
$$
CY_p=-\sum_{l=1}^{p-1}B(Y_l,Y_{p-l})
$$
holds.

Now we are ready to formulate sufficient conditions implying
existence of an implicit function which is determined by a
system of algebraic polynomial equations.

{\bf Theorem 1.} {\it Let $$ \sum_{p=0}^{q} Y_pt^p    \eqno(5)
$$ be an approximate solution of degree $q$ to a
system of algebraic polynomial equations $F(X)=0$.
Suppose that there exists a number $k$ $(0\leq k<q)$
such that, for all $i=1, 2, \dots ,q$ and $j=k, k+1,
\dots , q$, the equation $$ CY=-B(Y_i,Y_j)-B(Y_j,Y_i) $$ has
a solution which lies in the linear span of the vectors $Y_k$, $Y_{k+1}$,
\dots $Y_q$. Then $F(X)=0$ has an analytic family of solutions
$X(t)=\sum\limits_{p=0}^{\infty} X_pt^p$
whose initial coefficients coincide with the corresponding coefficients of
the approximate solution (5), i.e.,  $X_p=Y_p$ for each $p=0,1,\dots , q$.}

{\bf Proof.}
Denote by $L$ the linear span of the vectors
$Y_k$, $Y_{k+1}$, \dots ,$Y_q$.

Let $P$ be the set of all nonnegative integers $p$
for each of which there exists an approximate solution of degree $q+p$,
$$ \sum_{l=0}^{q+p} X_lt^l,\eqno(6) $$ to $F(X)=0$
such that (i) $X_l=Y_l$  for all $l=0, 1, \dots ,q$ and
(ii)  $X_l\in L$ for each $l=q+1, q+2, \dots , q+p$.

In view of the conditions of Theorem 1, $0\in P$.
Hence, $P\neq \emptyset$.
Verify that $P$ coincides with the set of all
nonnegative integers, $\mbox{\bf N}$.
It suffices to show that if $p\in P$ then $p+1\in P$.

So,  let $p\in P$  and let the approximate solution (6)
satisfy (i) and (ii).
To prove that $p+1\in P$, it is sufficient to find a vector
$X_{q+p+1}\in L$ satisfying the following system of linear
algebraic equations:
$$CX_{q+p+1}=-\sum_{l=1}^{q+p}B(X_l,X_{q+p+1-l}). \eqno(7) $$
According to (i) and (ii), each of the vectors $X_{q+1}, \dots
, X_{q+p}$ lies in $L$ and, hence,  each of the vectors $X_{q+1}, \dots
, X_{q+p}$ is a linear combination of the vectors
$Y_k$, $Y_{k+1}$, $\dots $, $Y_q$.
Therefore, the right-hand side of (7) is a linear combination of the vectors
$B(Y_i,Y_j)+B(Y_j,Y_i)$ with $1\leq i\leq q$ and $k\leq j\leq q$.
Then the conditions of Theorem 1 imply that there exists a solution to (7)
which lies in $L$.
Hence, $p+1\in P$ and $P=\mbox{\bf N}$.

Thus, we see that $F(X)=0$ has approximate solutions of arbitrarily high
degree whose initial coefficients coincide with the corresponding coefficients
of the approximate solution (5).
It remains to prove that these approximate solutions of
arbitrarily high degree can be used to construct an exact solution in the form of
a power series whose initial coefficients coincide with the corresponding
coefficients of (5).

Our approach is based on the following algebraic theorem by M.~Artin
(see \cite{Artin} and \cite{Sabitov-E}):
{\it Given a system of polynomial equations $f(x,y)=0$, where
$f=(f_1,\dots , f_k)$, $x=(x_1, \dots ,x_m)$, $y=(y_1,\dots ,y_n)$,
there exists an integer $\beta=\beta (m, n, d, \alpha)$
(depending on $m$, $n$, on the total degree $d$ of the polynomials $f$,
and on a nonnegative integer $\alpha$) such that if
$f(x,\overline{y}(x))=0 \ (\mbox{mod\,} x^\beta)$,
$\beta=\beta(m,n,d,\alpha)$, for some polynomial $\overline{y}(x)$
then the system $f(x,y)=0$ has a solution $y(x)$ that can be represented
in the form of a convergent power series whose coefficients
coincide with those of the polynomial $\overline{y}(x)$ up to the
term $x^\alpha$.}

To complete the proof of Theorem 1, we apply Artin's theorem
to $F(X)=0$ in the following way.
Put $\alpha = q$ and
let $\beta$ be an integer whose existence is provided by Artin's theorem.
From the above it follows that the approximate solution (5)
can be extended to an approximate solution of an arbitrarily high degree,
in particular, to an approximate solution of degree $\beta$.
Now Theorem 1 directly follows from Artin's theorem.

We now discuss several examples of using Theorem 1.

{\bf Example 1.}
Let $F:\mbox{\bf R}^3\to\mbox{\bf R}^3$ be given by (1), namely, let
$$
F_1(t,x_1,x_2,x_3)\equiv x_1^2+x_2^2-x_3^2-1 =0,
$$
$$
F_2(t,x_1,x_2,x_3)\equiv 3x_1+x_2-3x_3+1 =0,
$$
$$
F_3(t,x_1,x_2,x_3)\equiv x_1-3x_2+x_3+3 =0
$$
and let $X_0=(5,5,7)^{\rm T}$.
Direct calculations show that
$$
(\alpha_{ij}^{1})=
\left(
\begin{array}{ccc}
1 & 0 & 0 \\
0 & 1 & 0 \\
0 & 0 & -1
\end{array}
\right); \qquad
(\beta_{i}^{1})=
\left(
\begin{array}{c}
0 \\
0 \\
0
\end{array}
\right); \qquad
\gamma^1=-1;
$$
$$
(\alpha_{ij}^{2})=
\left(
\begin{array}{ccc}
0 & 0 & 0 \\
0 & 0 & 0 \\
0 & 0 & 0
\end{array}
\right); \qquad
(\beta_{i}^{2})=
\left(
\begin{array}{c}
3 \\
1 \\
-3
\end{array}
\right); \qquad
\gamma^2=1;
$$
$$
(\alpha_{ij}^{3})=
\left(
\begin{array}{ccc}
0 & 0 & 0 \\
0 & 0 & 0 \\
0 & 0 & 0
\end{array}
\right); \qquad
(\beta_{i}^{3})=
\left(
\begin{array}{c}
1 \\
-3 \\
1
\end{array}
\right); \qquad
\gamma^3=3;
$$
$$
B(X,Y)=(x_1y_1+x_2y_2-x_3y_3,0,0)^{\rm T};
$$
$$
A=\left(\begin{array}{ccc}
0 & 0 & 0 \\
3 & 1 & -3 \\
1 & -3 & 1
\end{array}\right); \qquad
C=\left(\begin{array}{ccc}
10 & 10 & -14 \\
3 & 1 & -3 \\
1 & -3 & 1
\end{array}\right); \qquad
\det C=0.
$$
Solving the homogeneous system of linear algebraic equations $CX=0$,
we find that the vector $X_1=(4,3,5)^{\rm T}$ constitutes a basis
for the space of its solutions.
Direct calculations show that $B(X_1,X_1)=(0,0,0)^{\rm T}$.
Hence, we can put $X_q=0$ for all $q\geq 2$.
Thus, Theorem 1 can be applied with $q=2$ and $k=1$.
It follows that $X_0$ is not an isolated solution to $F(X)=0$.
In the case under consideration, the corresponding family of solutions
may be explicitly written as $X(t)=X_0+tX_1$.
Its geometrical sense becomes obvious if we observe that the equation
$F_1(X)=0$ determines a one-sheet hyperboloid while the pair of linear
equations $F_2(X)=F_3(X)=0$ determines its straight line generator
which passes through the point $X_0$.

{\bf Example 2.}
Let $f:\mbox{\bf R}^2\to\mbox{\bf R}^1$ be given by the formula
$f(x_1,x_2)=x_1^3-x_2^2$ and let $X_0=(0,0)^{\rm T}$.

Transform the equation $f(X)=0$  to a system of equations each of which
is of degree 2:
$$
\begin{array}{l}
F_1(x_1,x_2,x_3)\equiv x_1x_3-x_2^2=0,\\
F_2(x_1,x_2,x_3)\equiv x_1^2-x_3=0.
\end{array}
\eqno(8)
$$
Then we have
$B(X,Y)=(\frac 12 x_1y_3-x_2y_2+\frac12 x_3y_1, x_1y_1)^{\rm T}$,
$$
C=
\left(
\begin{array}{ccc}
0 & 0 & 0 \\
0 & 0 & -1
\end{array}
\right).
$$
Instead of solving the corresponding equations (4) step by step,
we note that the equation $f(X)=0$ has an obvious analytic family
of solutions, $x_1=t^2$, $x_2=t^3$.
Whence we immediately obtain
$X_0=X_1=X_5=X_6=\dots =(0,0,0)^{\rm T}$,
$X_2=(1,0,0)^{\rm T}$, $X_3=(0,1,0)^{\rm T}$, $X_4=(0,0,1)^{\rm T}$.
Find the smallest values of $q$ and $k$
for which the hypotheses of Theorem 1 are fulfilled.

Direct calculations show that
$$
B(X_1,X_i)+B(X_i,X_1)=(0,0)^{\rm T}\qquad\mbox{for all $i\geq 1$,}
$$
$$
B(X_2,X_i)+B(X_i,X_2)=(0,0)^{\rm T}\qquad\mbox{for $i=3$ or $i\geq 5$,}
$$
$$
B(X_3,X_i)+B(X_i,X_3)=(0,0)^{\rm T}\qquad\mbox{for all $i\geq 4$,}
$$
$$
B(X_4,X_i)+B(X_i,X_4)=(0,0)^{\rm T}\qquad\mbox{for all $i\geq 4$,}
$$
$$
B(X_2,X_2)=(0,1)^{\rm T},\quad
B(X_3,X_3)=(-1,0)^{\rm T},
$$
$$
B(X_2,X_4)+B(X_4,X_2)=(1,0)^{\rm T}.
$$
Hence, the hypotheses of Theorem 1 are fulfilled for $q=k=5$
and are not fulfilled for any smaller values of $q$ and $k$.
Thus, as soon as we obtain the approximate solution
$X_0+tX_1+t^2X_2+t^3X_3+t^4X_4+t^5X_5$,
we can assert that it can be extended to an exact solution to (8).
However, we cannot make the same assertion by using only the approximate
solution $X_0+tX_1+t^2X_2+t^3X_3+t^4X_4$.

Now we give an example of a system of algebraic equations
possessing an analytic family of solutions that cannot
be obtained by Theorem 1 for any values of $q$ and $k$.

{\bf Example 3.}
Let $F:\mbox{\bf R}^3\to\mbox{\bf R}^2$ be given by the formulas
$$
\begin{array}{l}
F_1(x_1,x_2,x_3)\equiv x_1^2+x_2^2+x_3^2-4,\\
F_2(x_1,x_2,x_3)\equiv (x_1-1)^2+x_2^2-1
\end{array}
\eqno(9)
$$
and let $X_0=(2,0,0)^{\rm T}$.
Direct calculations show that
$$
(\alpha_{ij}^{1})=
\left(
\begin{array}{ccc}
1 & 0 & 0 \\
0 & 1 & 0 \\
0 & 0 & 1
\end{array}
\right); \qquad
(\beta_{i}^{1})=
\left(
\begin{array}{c}
0 \\
0 \\
0
\end{array}
\right); \qquad
\gamma^1=-4;
$$
$$
(\alpha_{ij}^{2})=
\left(
\begin{array}{ccc}
1 & 0 & 0 \\
0 & 1 & 0 \\
0 & 0 & 0
\end{array}
\right); \qquad
(\beta_{i}^{2})=
\left(
\begin{array}{c}
-2 \\
0 \\
0
\end{array}
\right); \qquad
\gamma^2=0;
$$
$$
B(X,Y)=(x_1y_1+x_2y_2+x_3y_3,x_1y_1+x_2y_2)^{\rm T};
$$
$$
A=\left(\begin{array}{ccc}
0 & 0 & 0 \\
-2 & 0 & 0
\end{array}\right); \qquad
C=\left(\begin{array}{ccc}
4 & 0 & 0 \\
2 & 0 & 0
\end{array}\right).
$$

Note that $\mbox{rank\,} C = 1$,
$\mbox{im\,} C = \{ (\xi , \eta )\in \mbox{\bf R}^2 \vert \xi = 2\eta \}$,
$\mbox{ker\,} C = \{ (u,v,w)\in\mbox{\bf R}^3 \vert u = 0\}$, and
$\mbox{dim ker\,} C = 2$.

Obviously, (9) determines the Viviani curve and thus
admits the following analytic family of solutions:
$$
\begin{array}{l}
x_1(t)=1+\cos t, \\
x_2(t)= \sin t, \\
x_3(t)=2\sin(t/2),\\
X(t)=(x_1(t),x_2(t),x_3(t))=\sum\limits_{p=0}^{\infty}t^pX_p.
\end{array}
$$

It is clear that each expression $\sum\limits_{p=0}^{N}t^pX_p$
is an approximate solution of some degree to (9).
Suppose it satisfies the hypotheses of Theorem 1 with some $q$ and $k$.

From the proof of Theorem 1 it follows that, for $p\geq k$, the vector
$X_p$ is constructed as a linear combination of solutions to the equations $$
CX=-B(X_i,X_j)-B(X_j,X_i) $$ with $1\leq i\leq q$ and $k\leq j\leq
q$. Hence, each vector $B(X_2,X_p)+B(X_p,X_2)$
must lie in the image of the operator $C$, i.e., its first coordinate
must be twice as large as its second coordinate.
However, this condition is not satisfied, since
$B(X_2,Y)+B(Y,X_2)=(-y_1,-y_1)^{\rm T}$ and infinitely many $X_p$
have nonzero first coordinate.

Thus, example 3 shows that the conditions of Theorem 1 are not necessary
for existence of an implicit function.
This means that Theorem 1 cannot be used for proving that a given
solution to a system of algebraic equations is isolated.
In the next section we give several additional conditions under which
the conditions of Theorem 1 are not only sufficient but also necessary
for existence of an analytic family of solutions to a system of
algebraic equations.

\section{Necessary conditions for existence of an implicit function}

A primary necessary condition for existence of a continuous family
of solutions is known in the theory of bending of smooth surfaces
at least since S.~Cohn-Vossen \cite{Cohn-Vossen}.
For systems of algebraic equations it can be formulated as follows
(we use the notation introduced in the previous section):

{\bf Theorem 2.} {\it If zero is the only solution to the system $CX=0$
then the system of algebraic equations $F(X)=0$ has no nonconstant
analytic family of solutions representable in the form of a convergent power
series with prescribed initial term $X_0$.}

{\bf Proof} is carried out by way of contradiction. Suppose that $F(X)=0$
has a nonconstant analytic family of solutions which is represented
in the form of a convergent power series
$$X(t)=\sum\limits_{p=0}^{\infty}t^pX_p $$ and suppose that $q$
is the smallest positive number such that $X_q\neq 0$.
According to (4), $X_q$ satisfies the following system of linear
algebraic equations: $$CX_q=-\sum_{p=1}^{q-1}B(X_p,X_{q-p})=0.$$
However, due to the hypotheses of Theorem 2, the last system
has only zero solution. Hence, $X_q=0$.
This contradiction proves Theorem 2.

In the theory of bending of smooth surfaces some more advanced
necessary conditions are also known (see, for example, \cite{Con-Adv},
\cite{Efimov-52}).
An algebraic version is stated in the following theorem:

{\bf Theorem 3.} {\it If the system of algebraic equations
$F(X)=0$ and the vector $X_0$ are such that no approximate solution
of the first degree $X_0+tX_1$, with $X_1\neq 0$, can be extended
to an approximate solution of the second degree, then
$F(X)=0$ has no nonconstant analytic family of solutions representable in
the form of a convergent power series with initial term $X_0$.}

{\bf Proof} is carried out by way of contradiction.
Suppose that $F(X)=0$
has a nonconstant analytic family of solutions which is represented
in the form of a convergent power series
$$X(t)=\sum\limits_{p=0}^{\infty}t^pX_p $$ and suppose that $q$
is the smallest positive number such that $X_q\neq 0$.
Then $X_q$ lies in the kernel of $C$ and, according to the hypotheses
of Theorem 3, it follows that $B(X_q,X_q)$ does not lie in the kernel of $C$.
In view of (4), $X_{2q}$ satisfies the following system of linear equations:
$$CX_{2q}=-\sum_{p=1}^{2q-1}B(X_p,X_{2q-p})=-B(X_q,X_q).$$
Since $B(X_q,X_q)$ does not lie in the kernel of $C$,
the last system has no solutions.
This contradiction proves Theorem 3.

We now discuss several additional conditions under which the
conditions of Theorem 1 are not only sufficient but also necessary
for existence of an analytic family of solutions to a system of
algebraic equations.
The results we aim to obtain will generalize Theorems 2 and 3 in a sense.
We start with the case in which there are few linearly independent vectors
in the sequence $X_1, X_2, \dots , X_q$ of coefficients of
approximate solutions.

{\bf Theorem 4.} {\it
Suppose that the system of algebraic equations $F(X)=0$
has an analytic family of solutions which is represented
in the form of a convergent power series
$$X(t)=\sum\limits_{p=0}^{\infty}t^pX_p $$ and suppose
that the vectors $X_3$ and $X_4$ belong to the linear span
of the vectors $X_1$ and $X_2$.
Then, for all $1\leq i,j\leq 2$, the system of linear equations $$
CX=-B(X_i,X_j)-B(X_j,X_i) $$ has a solution that belongs to the
linear span of the vectors $X_1$ and $X_2$.}

{\bf Proof.}
The vectors $X_1$, $X_2$, $X_3$, and $X_4$
satisfy the following equations:
$$
\begin{array}{l}
CX_1=0,\\
CX_2=-B(X_1,X_1),\\
CX_3=-B(X_1,X_2)-B(X_2,X_1),\\
CX_4=-B(X_1,X_3)-B(X_2,X_2)-B(X_3,X_1).
\end{array}
$$
Let $L$ denote the linear span of $X_1$ and $X_2$.
The second equation implies $B(X_1,X_1)\in CL$.
Since $X_3\in L$, the third equation implies $B(X_1,X_2)+B(X_2,X_1)\in CL$.
On the other hand, $X_3\in L$ implies that  $X_3=c_3^1X_1+c_3^2X_2$,
with some reals $c_3^1$ and $c_3^2$.
Hence, the fourth equation can be written as
$CX_4=-2c_3^1B(X_1,X_1)-c_3^2[B(X_1,X_2)+B(X_2,X_1)]-B(X_2,X_2)$.
Here $CX_4$ belongs to $CL$ according to the hypotheses of Theorem 4,
$B(X_1,X_1)$ and $B(X_1,X_2)+B(X_2,X_1)$ belong to $CL$
in view of the above proof.
Therefore, $B(X_2,X_2)\in CL$. This completes the proof of Theorem 4.

{\bf Theorem 5.} {\it
Suppose that the system of algebraic equations $F(X)=0$
has an analytic family of solutions which is represented
in the form of a convergent power series
$$X(t)=\sum\limits_{p=0}^{\infty}t^pX_p $$ and suppose that the
vectors $X_4$, $X_5$, $X_6$, and $X_7$ belong to the linear span
of the vectors $X_1$, $X_2$, and $X_3$.
Then, for all $1\leq i,j\leq 3$, the system of linear equations $$
CX=-B(X_i,X_j)-B(X_j,X_i) $$ has a solution that belongs to the
linear span of the vectors $X_1$, $X_2$, and $X_3$.}

{\bf Proof.} Let $\alpha\in\mbox{\bf R}$ be an arbitrary real.
Change the variable $t=\tau +\alpha\tau^2$ in the solution
$X(t)=\sum\limits_{p=0}^{\infty}t^pX_p$ : $Y(\tau)\equiv X(\tau
+\alpha\tau^2)=\sum\limits_{p=0}^{\infty}\tau^pY_p$.
Clearly, $Y(\tau )$ is an analytic family of solutions to the equation
$F(Y)=0$ and, thus, the equality
$$CY_q=-\sum_{p=1}^{q-1}B(Y_p,Y_{q-p})\eqno(10) $$
holds for each $q\geq 1$.

On the other hand, $Y_p$ can be written in terms of $X_i$
by collecting similar terms in the expression
$$
\sum\limits_{p=0}^{\infty}\tau^pY_p=
\sum\limits_{p=0}^{\infty}(\tau+\alpha\tau^2)^pY_p.
$$
We thus obtain
$$
Y_0=X_0;
$$
$$
Y_1=X_1;\eqno(11)
$$
$$
Y_2=X_2+\alpha X_1;\eqno(12)
$$
$$
Y_3=X_3+2\alpha X_2; \eqno(13)
$$
$$
Y_4=X_4+3\alpha X_3 +\alpha^2 X_2;\eqno(14)
$$
$$
Y_5=X_5+4\alpha X_4 +3\alpha^2 X_3;\eqno(15)
$$
$$
Y_6=X_6+5\alpha X_5 +6\alpha^2 X_4 +\alpha^3 X_3.\eqno(16)
$$

According to the hypotheses of Theorem 5, each of the vectors
$X_4$, $X_5$, and $X_6$ belongs to the linear span of the vectors
$X_1$, $X_2$, and $X_3$.
It follows that   $X_j=c_j^1X_1+c_j^2X_2+c_j^3X_3$ for each
$4\leq j\leq 6$, with some reals $c_j^i$
($1\leq i\leq 3$, $4\leq j\leq 6$).
Taking  these formulas into account, we can write (14)--(16)
in the following form:
$$
\begin{array}{l}
Y_4=(c_4^3+3\alpha)X_3+(c_4^2+\alpha^2)X_2+c_4^1X_1,\\
Y_5=(c_5^3+4\alpha c_4^3 +3\alpha^2)X_3+(c_5^2+4\alpha c_4^2)X_2
+(c_5^1+4\alpha c_4^1)X_1,\\
Y_6=(c_6^3+5\alpha c_5^3 +6\alpha^2 c_4^3+\alpha^3)X_3+
(c_6^2+5\alpha c_5^2+6\alpha^2c_4^2)X_2\\
\phantom{Y_6=}+(c_6^1+5\alpha c_5^1 +6\alpha^2 c_4^1)X_1.
\end{array}\eqno(17)
$$

Let $L$ denote the linear span of $X_1$, $X_2$, and $X_3$.

For the case $q=2$ equation (10) gives $CY_2=-B(Y_1,Y_1)$.
Taking into account (11) and (12), we conclude from here that
$CX_2+\alpha CX_1=-B(X_1,X_1)$, and finally that $B(X_1,X_1)\in CL$.

For the case $q=3$ equation (10) gives $CY_3=-B(Y_1,Y_2)-B(Y_2,Y_1)$.
Taking into account (12) and (13), we conclude from here that
$CX_3+2\alpha CX_2=-2\alpha B(X_1,X_1)- [B(X_1,X_2)+B(X_2,X_1)]$,
and finally that $B(X_1,X_2)+B(X_2,X_1)\in CL$.

For the case $q=4$, equation (10) yields
$CY_4=-B(Y_1,Y_3)-B(Y_2,Y_2)-B(Y_3,Y_1)$.
Taking (14) and (17) into account, we conclude that
$(c_4^3+3\alpha)CX_3+(c_4^2+\alpha^2) CX_2 + c_4^1 CX_1 =
-\alpha^2 B(X_1,X_1)-3\alpha [B(X_1,X_2)+B(X_2,X_1)]
-[B(X_1,X_3)+B(X_3,X_1)]-B(X_2,X_2)$
and, finally,
$$
[B(X_1,X_2)+B(X_2,X_1)]+B(X_2,X_2)\in CL.\eqno(18)
$$

Similarly, from (10), we obtain in the case $q=5$,
$$
[B(X_1,X_3)+B(X_3,X_1)]+[B(X_2,X_3)+B(X_3,X_2)]\in CL, \eqno(19)
$$
in the case $q=6$,
$$
\begin{array}{l}
(c_5^3+4\alpha c_4^3 +3\alpha^2)[B(X_1,X_3)+B(X_3,X_1)]\\
+(2c_4^2-10\alpha^2-4\alpha c_4^3)B(X_2,X_2)\\
+(c_4^3+3\alpha)[B(X_2,X_3)+B(X_3,X_2)]\in CL,
\end{array}\eqno(20)
$$
and, in the case $q=7$,
$$
\begin{array}{l}
(c_4^1-\alpha c_4^2 +6\alpha^4+c_6^3+5\alpha c_5^3+8\alpha^2 c_4^3)
[B(X_1,X_3)+B(X_3,X_1)]\\
+(c_5^2+4\alpha c_4^2-2\alpha c_5^3+8\alpha^2 c_4^3-6\alpha^3)B(X_2,X_2)\\
+(c_5^3+2\alpha c_4^3 +c_4^2 + 4\alpha^2)[B(X_2,X_3)+B(X_3,X_2)]\\
+(c_4^3+3\alpha)B(X_3,X_3)\in CL.
\end{array}\eqno(21)
$$

We may treat relations (18)--(21) as a system of algebraic linear equations
with respect to the following four vector-valued variables:
$B(X_1,X_3)+B(X_3,X_1)$, $B(X_2,X_2)$, $B(X_2,X_3)+B(X_3,X_2)$,
and $B(X_3,X_3)$.
The right-hand sides of the corresponding equations are some vectors in $CL$.
Denote them by $U_1$, $U_2$, $U_3$, and $U_4$.
It suffices to prove that the determinant of the system is not equal to zero:
in this case, each of the four vector-valued variables may be represented
as a linear combination of  $U_1$, $U_2$, $U_3$, and $U_4$ and, thus,
lies in $CL$.

Direct calculations show that the determinant of the system corresponding
to (18)--(21) equals
$$
\begin{array}{l}
-6\alpha^4 +18\alpha^3 +20\alpha^2 c_4^3 +
\alpha [2c_4^2 +(c_4^3)^2 -c_5^3]\\
+[-c_4^1 +c_5^2 -c_4^2 c_4^3 -(c_4^3)^2 +2 c_4^3 c_5^3 -c_6^3].
\end{array}
$$
Obviously, this polynomial in $\alpha$ cannot be equal to zero identically
for any values of the coefficients $c_j^i$.
Hence, we can find a value of $\alpha$ such that
the determinant of the system corresponding to (18)--(21) do not vanish.
This completes the proof of Theorem 5.

The above example 3 shows that, if the four vectors $X_1$, $X_2$, $X_3$, and
$X_4$ are linearly independent, it may occur that some of the vectors
$B(X_i,X_j)+B(X_j,X_i)$ do not lie in the image of $C$ while the system
$F(X)=0$ determines some nonconstant implicit function.
This means that there is no direct generalization of Theorems 4 and 5
to the case in which the four vectors $X_1$, $X_2$, $X_3$, and
$X_4$ are linearly independent.

Obviously, this circumstance hampers proving that
a system $F(X)=0$ admits no analytic family of solutions.
There are also some other obstacles in proving this.
The first obstacle is of purely technical nature and consists in
a considerable increase of calculations: if we find that a first-order
approximate solution  $X_0+tX_1$ admits an extension  $X_0+tX_1+t^2X_2$
to some second-order approximate solution then,
for every $\widetilde{X}\in\mbox{ker\, } C$, the expression
$X_0+tX_1+t^2(X_2+\widetilde{X})$ also gives us some second-order
approximate solution; hence, we need to study the possibility
of extending a second-order approximate solution to a third-order
approximate solution for some family of second-order solutions
rather than for a single solution.
Another reason is of more principle character.
Suppose that we are able to make our way through the above-described
increase of calculations and suppose we find a number $N$
such that no first-order approximate solution
$X_0+tX_1$, $X_1\in\mbox{ker\, } C$, $X_1\neq 0$,
can be extended to any approximate solution of order $N$.
Can we conclude that the system $F(X)=0$ defines no implicit
function in a neighborhood of the point $X_0$? No, we cannot!
We have to verify that there is no such an extension
either for   $X_0+t\cdot 0+t^2X_1$, or for $X_0+t\cdot 0+t^2\cdot 0
+t^3 X_1$, or for any other approximate solution which has
several vanishing initial coefficients
(here $X_1\in\mbox{ker\, } C$, $X_1\neq 0$).
Example 2 shows that a system $F(X)=0$ may admit an exact solution
with several vanishing initial coefficients $X_p$.
On the other hand, we have no estimation for the number of zero
coefficients in the Maclaurin expansion of an implicit function
defined by a system $F(X)=0$.
Hence, we have to study an infinite set of cases caused by
``writing zeros'' at initial positions of an approximate solution.
Therefore, in general, we have no algorithm which can guarantee
absence of an implicit function.

Nevertheless, below we will show that such an algorithm does exist in
the case $\mbox{dim~ker\, } C=1$.
First of all, we make the terminology more accurate.
As before, let $F(X)=0$ be a system of algebraic equations
each of which is of degree 1 or 2.
Let a bilinear operator $B$ and a linear operator $C$ be constructed
by means of the system.
Suppose that $\mbox{dim~ker\, } C=1$.
In the domain of $C$ fix a codimension-1 subspace $T$ such that
$T\cap\mbox{ker\, } C=\{ 0\}$.
A formal power series $X(t)=\sum\limits_{p=0}^{\infty}t^pX_p$
is said to be a {\it $T$-standard formal solution}
to $F(X)=0$ if the following conditions are fulfilled:

1) $CX_q=-\sum\limits_{p=1}^{q-1}B(X_p,X_{q-p})$ for each  $q\geq 1$;

2) $X_1\neq 0$;

3) $X_p\in T$ for each $p\geq 2$.

The following theorem plays a key role in our approach:

{\bf Theorem 6.} {\it Let the system $F(X)=0$ admit an exact
nonconstant solution which can be represented in the form of
a convergent power series
$X(t)=\sum\limits_{p=0}^{\infty}t^pX_p$, let $\mbox{dim~ker\, } C=1$,
and let $T$ be a codimension-1 subspace such that
$T\cap\mbox{ker\, } C=\{ 0\}$.
Then $F(X)=0$ admits a $T$-standard formal solution
$Y(t)=\sum\limits_{p=0}^{\infty}t^pY_p$ such that $Y_0=X_0$.}

To avoid interrupting our presentation, we prove Theorem 6
at the end of this section.

The coefficients $Y_p$ of a $T$-standard formal solution
to $F(X)=0$ can be found as solutions to the following
system of linear algebraic equations:
$CY_p=-\sum\limits_{l=1}^{p-1}B(Y_l,Y_{p-l})$.
The condition  $Y_p\in T$ implies that the solution $Y_p$  is unique
(if existent, of course) and thus no increase of calculations occurs.
On the other hand, if, for some $p$, a solution $Y_p$ does not exist
then there is no $T$-standard formal solution to $F(X)=0$.
According to Theorem 6, this implies that $F(X)=0$
has no exact nonconstant solution in the form of a convergent
power series (with arbitrarily many vanishing initial
coefficients).
Thus, we have a finite algorithm which, in some cases,
can guarantee absence of an implicit function determined by $F(X)=0$
in a neighborhood of the point $X_0$.
We present a test example of executing the algorithm proposed.

{\bf Example 4.}
Let $F:\mbox{\bf R}^3\to\mbox{\bf R}^3$ be given by the formulas
$$
\begin{array}{l}
F_1(x_1,x_2,x_3)\equiv x_1^2+x_2^2+x_3^2-4,\\
F_2(x_1,x_2,x_3)\equiv (x_1-3)^2+x_2^2-1,\\
F_3(x_1,x_2,x_3)\equiv x_2
\end{array}
$$ and let $X_0=(2,0,0)^{\rm T}$.
It is clear that the equation $F_1=0$
defines a sphere in $\mbox{\bf R}^3$
while $F_2=0$ defines a cylinder which has a single
common point with the sphere, namely, the point $X_0$.
Thus $F(X)=0$ defines no implicit function.
Demonstrate how Theorem 6 can be used to reach the same conclusion.

Direct calculations show that
$$
(\alpha_{ij}^{1})=
\left(
\begin{array}{ccc}
1 & 0 & 0 \\
0 & 1 & 0 \\
0 & 0 & 1
\end{array}
\right); \qquad
(\beta_{i}^{1})=
\left(
\begin{array}{c}
0 \\
0 \\
0
\end{array}
\right); \qquad
\gamma^1=-4;
$$
$$
(\alpha_{ij}^{2})=
\left(
\begin{array}{ccc}
1 & 0 & 0 \\
0 & 1 & 0 \\
0 & 0 & 0
\end{array}
\right); \qquad
(\beta_{i}^{2})=
\left(
\begin{array}{c}
-6 \\
0 \\
0
\end{array}
\right); \qquad
\gamma^2=-1;
$$
$$
(\alpha_{ij}^{3})=
\left(
\begin{array}{ccc}
0 & 0 & 0 \\
0 & 0 & 0 \\
0 & 0 & 0
\end{array}
\right); \qquad
(\beta_{i}^{3})=
\left(
\begin{array}{c}
0 \\
1 \\
0
\end{array}
\right); \qquad
\gamma^3=0;
$$
$$
B(X,Y)=(x_1y_1+x_2y_2+x_3y_3,x_1y_1+x_2y_2,0)^{\rm T};
$$
$$
A=\left(\begin{array}{ccc}
0 & 0 & 0 \\
-6 & 0 & 0\\
0 & 1 & 0
\end{array}\right); \qquad
C=\left(\begin{array}{ccc}
4 & 0 & 0 \\
-2 & 0 & 0 \\
0 & 1 & 0
\end{array}\right).
$$

Note that $\mbox{rank\,} C = 2$,
$\mbox{im\,} C = \{ (\xi , \eta , \zeta)\in \mbox{\bf R}^3 \vert \xi = 2\eta \}$,
$\mbox{ker\,} C = \{ (u,v,w)\in\mbox{\bf R}^3 \vert u = v = 0\}$,
$\mbox{dim ker\,} C = 1$,
$X_1=(0,0,1)$, and
$B(X_1,X_1)=(1,0,0)^{\rm T}\notin \mbox{im\, } C$.

The latter relation implies that the approximate solution
$X_0+tX_1$ cannot be extended to any approximate solution
of the second order.
Hence, there is no $T$-standard formal solution to $F(X)=0$.
By Theorem 6, $F(X)=0$ defines no implicit function
in a neighborhood of $X_0$.

{\bf Proof of Theorem 6.} Let $N$ be the least positive
integer $p$ such that $X_p\neq 0$.
Let $q$ be the greatest integer such that the exact nonconstant
solution to $F(X)=0$ (which, according to the claims of Theorem 6,
is representable in the form of a convergent power series
$X(t)=\sum\limits_{p=0}^{\infty}t^pX_p$) possesses the following
properties:
(i) $X_p=0$ for all $0<p\leq q$, $p\neq 0\ (\mbox{mod\, } N)$ and
(ii) $X_p\in T$ for all $0<p\leq q$, $p\neq N$, $p=0\ (\mbox{mod\, } N)$.

We will verify that there exists a polynomial change of variables
$t=t(\tau)$ such that the new exact nonconstant solution
$\widetilde{X}(\tau) \equiv X(t(\tau))$ to $F(X)=0$,
which is representable in the form of a convergent power series
$\widetilde{X}(\tau)=\sum\limits_{p=0}^{\infty}\tau^p\widetilde{X}_p$,
possesses properties (i) and (ii) for $q+1$ and is such that
$\widetilde{X}_0=X_0$, $\widetilde{X}_N=X_N$.

In other words, we will prove that, by means of polynomial
changes of variable $t$, the coefficients $X_p$ can be transformed
one after another into the zero vector if $p$ is not divisible by $N$
and into some vectors lying in $T$ if $p$ is divisible by $N$
in such a manner that  $X_0$ and $X_N$  remain unchanged and,
after such a transformation, we again obtain an exact solution to
$F(X)=0$ representable in the form of a convergent power series.
Accomplishing infinitely many such polynomial changes of variable $t$
and not controlling the radii of convergence of power series which appear
in this process, we obtain a formal power series
$Z(\tau)=\sum\limits_{p=0}^{\infty}\tau^pZ_p$
whose coefficients $Z_p$ ($p=0,1,\dots $) possess the following properties:

(a) $Z_0=X_0$;

(b) $Z_p=0 \mbox{\ for all\ } p\neq 0 \ (\mbox{mod\, } N)$;

(c) $Z_N=X_N\neq 0$;

(d) $Z_p\in T \mbox{\ for all\ } p>0,\ p= 0 \ (\mbox{mod\, } N)$;

(e) $CZ_q=-\sum\limits_{p=1}^{q-1}B(Z_p,Z_{q-p})$
for every $q\geq 1$.

Finally, executing the change of variable  $t=\tau^N$ in the
formal power series $Z(\tau)$, we obtain  a $T$-standard formal
solution $Y(t)$ whose existence is asserted in Theorem 6.

So, to complete the proof of Theorem 6, it remains to prove
that a solution $X(t)$ possessing
properties (i) and (ii) for some $q$, can be transformed into a solution
$\widetilde{X}(\tau)$ possessing properties (i) and (ii)
for $q+1$ and satisfying $\widetilde{X}_0=X_0$, $\widetilde{X}_N=X_N$.

Put $q=iN+j$, where $0\leq j\leq N-1$.
Consider a change of variable $t=\tau +\alpha\tau^{q+1-N}$
(here $\tau$ is a new variable and $\alpha$ is a constant
whose value will be specified later): $$
\begin{array}{ll}
\widetilde{X}(\tau) &= X(\tau+\alpha\tau^{q+1-N})\\
 &=X_0 + (\tau + \alpha \tau^{q+1-N})^NX_N +
(\tau + \alpha \tau^{q+1-N})^{2N}X_{2N} + \cdots \\
 & \hspace{9mm}+ (\tau + \alpha \tau^{q+1-N})^{iN}X_{iN} +(\tau + \alpha
\tau^{q+1-N})^{iN+j}X_{iN+j}+ \cdots \\
 &= X_0+\tau^NX_N+\tau^{2N}X_{2N} +\cdots +\tau^{iN}X_{iN}+
 \tau^q(X_q+N\alpha X_N)+\cdots .
\end{array}$$

First, consider the case $0< j\leq N-1$.
We know that $X_p$ is a solution to the following system of linear
algebraic equations:
$$ CX_q=-\sum_{r=1}^{iN+j}B(X_r,X_{iN+j-r}). \eqno(22) $$
If $r= 0\ (\mbox{mod\, } N)$ then $iN+j-r =j \ (\mbox{mod\, } N)$
and, in particular, $iN+j-r\neq 0\ (\mbox{mod\, } N)$.
Hence, for arbitrary $1\leq r\leq iN+j-1$, either $X_r=0$ or
$X_{iN+j-r}=0$.
Therefore, the right-hand side of (22) equals zero and $X_q\in\mbox{ker\,} C$.
On the other hand, $X_N\in\mbox{ker\,} C$ and $\mbox{dim ker\,} C = 1$.
Consequently, the vectors $X_q$ and $X_N$ are collinear.
Since $X_N\neq 0$, there exists $\alpha$ such that $X_q+\alpha N X_N=0$.
Under such a choice of $\alpha$,
the exact solution $\widetilde{X}(\tau)=X(\tau+\alpha\tau^{q+1-N})$
possesses properties (i) and (ii) for $q+1$
as well as $\widetilde{X}_0=X_0$ and $\widetilde{X}_N=X_N$.

Now consider the case $j=0$.
In this case, if $r=sN$ $(0\leq s\leq i)$ then $iN+j-r=(i-s)N$.
Hence, (22) can be rewritten as
$$ CX_q=-\sum_{s=1}^{i}B(X_{sN},X_{(i-s)N}). $$
Generally speaking, the right-hand side of the last expression does not
equal zero.
So, in general, $X_q$ does not lie in $\mbox{ker\,} C$.
Nevertheless, using the fact that the linear span of the subspaces
$T$ and $\mbox{ker\,} C$ coincides with the whole space,
we can find $\alpha$ such that $X_q+\alpha N X_N\in T$.
Under such a choice of $\alpha$,
the exact solution $\widetilde{X}(\tau)=X(\tau+\alpha\tau^{q+1-N})$
possesses properties (i) and (ii) for $q+1$
as well as $\widetilde{X}_0=X_0$ and $\widetilde{X}_N=X_N$.

Thus, the possibility has been proven of transferring the solution
$X(t)$ into a solution $\widetilde{X}(\tau)$
which possesses properties (i) and (ii) as well as
$\widetilde{X_0}=X_0$, $\widetilde{X}_N=X_N$.

This completes the proof of Theorem 6.

\section{Applications to studying flexible polyhedra and frameworks}

Let $K$ be a simplicial complex whose body is an $(n-1)$-dimensional
connected compact topological manifold without boundary.
A {\it polyhedron} in the $n$-dimensional Euclidean space
$\mbox{\bf R}^n$ is, by definition, a continuous mapping
$f:K\to \mbox{\bf R}^n$ which is linear on each simplex.
Sometimes, the image of $K$ under $f$ is also referred to as a polyhedron.
By a {\it polyhedral sphere} in  $\mbox{\bf R}^n$ we mean a polyhedron
$f:K\to \mbox{\bf R}^n$ with the body of $K$ homomorphic to the sphere.

We say that a polyhedron has no self-intersections if the mapping $f$
is (globally) injective.
In the present article, we consider polyhedra both with and without
self-intersections.

{\bf Definition.} A polyhedron $P=f(K)$ is {\it flexible} if there
exists a family of polyhedra $P_t=(f_t,K)$, $0\leq t \leq 1$,
which is analytic with respect to the parameter $t$ and for which
the following conditions are satisfied:

1) $P=P_0$;

2) for arbitrary 0-dimensional simplices $v_j$ and $v_k$ of $K$
belonging to a 1-dimensional simplex of $K$, the equality
$|f(v_j)-f(v_k)|=|f_t(v_j)-f_t(v_k)|$ holds for all $0\leq
t\leq 1$ (henceforth $|y|$  stands for
the Euclidean norm of a vector
$y=(y_1,y_2,\dots ,y_n)\in \mbox{\bf R}^n$, i.e., $|y|^2=y_1^2+
y_2^2 +\dots +y_n^2)$;

3) there exist two 0-dimensional simplices $v_j$ and $v_k$ of $K$
which do not belong to any 1-dimensional simplex of $K$ and
for which the expression $|f_t(v_j)-f_t(v_k)|$ is not constant
in $t\in [0,1]$.

A family $P_t$ with the above properties 1)--3) is a {\it
nontrivial flexion} of $P$.
Note that the simplicial complex $K$ remains constant during the flexion.

In other words, a polyhedron is flexible if its spatial shape
can be changed analytically with respect to a parameter (see
condition 3)) without changing its intrinsic metrics (see condition 2)).
Incidentally, the analyticity requirement with respect to the parameter
can be substantially weekend.
Namely, in \cite{Gluck} it is shown that if there exists
a deformation of a polyhedron which is continuous with respect to a
parameter and obeys the above conditions 1)--3) then there also exists
a deformation of the polyhedron which is analytic with respect to (possibly)
another parameter and obeys conditions 1)--3).

During the last 25 years, the following two remarkable results were
obtained in the theory of flexible polyhedra:
in 1977 R.~Connelly gave an example of a flexible polyhedral
sphere in $\mbox{\bf R}^3$ without self-intersections \cite{Connelly},
and in 1996 I.~Kh.~Sabitov published a complete proof of the statement
asserting that each flexible polyhedron (even with self-intersections)
in $\mbox{\bf R}^3$ preserves the (oriented) volume bounded by it
during the flexion \cite{Sabitov-FPM}.
The last statement was known during several decades as the
``Bellows conjecture''.
Other proofs of this statement may be found in
\cite{Con-Sab-Walz}, \cite{Sabitov-Sb}, and \cite{Sabitov-DCG}.

We are interested in the problem of ``practical'' recognition whether
a given polyhedron $f:K\to \mbox{\bf R}^n$ is flexible or not.
It is sufficient to describe a flexion  $f_t:K\to \mbox{\bf R}^n$
of the polyhedron via motions of its 0-dimensional simplices
$v_j$ $(1\leq j\leq N)$: $x_j(t)=f_t(v_j)\in \mbox{\bf R}^n$.
Furthermore, the above conditions 1)--3) can be reformulated as follows:
1') the vectors $x_j(0)$ are given;
2') if 0-dimensional simplices $v_i$ and $v_j$ are joint together in $K$
by a 1-dimensional simplex then the equality
$|x_i(t)-x_j(t)|^2=|x_i(0)-x_j(0)|^2$ holds true for all
$0\leq t\leq 1$;
3') there exist 0-dimensional simplices $v_i$ and $v_j$
that are not joint together in $K$ by a 1-dimensional simplex
and such that the expression $|x_i(t)-x_j(t)|^2$
is not constant in $t\in [0, 1]$.

In other words, the problem of whether a given polyhedron is flexible or not
is equivalent to the problem of whether the set of vectors
$x_j=x_j(0)$, $j=1,2,\dots ,N$, is an isolated solution to
the system of algebraic equations
$$|x_i-x_j|^2=|x_i(0)-x_j(0)|^2\eqno(23)$$
or this system defines an implicit function $x_j=x_j(t)$
in a neighborhood of the point $x_j=x_j(0)$, $j=1,2,\dots ,N$.
It is worth bearing in mind that we are not interested in motions of
$f(K)$ as a rigid body in $\mbox{\bf R}^n$, i.e.,
we are only looking for solutions with properties 3) or 3').
This requirement can be easily satisfied in the following manner.
Fix an $(n-1)$-dimensional simplex of $K$.
Suppose its 0-dimensional simplices are denoted by
$v_1$, $v_2$, \dots , $v_n$.
We agree that $f_t(v_1)$ lies at the origin during the course of
deformation of the polyhedron (and, thus, always has zero coordinates);
$f_t(v_2)$ always lies on the first coordinate axis in $\mbox{\bf R}^n$
(and, thus, all but the first coordinates of it vanish identically);
$f_t(v_3)$ always lies in the 2-dimensional plane spanned by the first and
second coordinate axis in $\mbox{\bf R}^n$
(and, thus, all but the first and second coordinates of it vanish
identically); and so on.
This construction reduces the number of independent variables in (23)
but does not preclude the above reduction of the decision problem of
whether the polyhedron is flexible to the problem of whether the given
solution to the system of algebraic equations is isolated or
the system defines an implicit function in a neighborhood of the solution.

The decision problem of whether a framework is flexible in $\mbox{\bf R}^n$
can be reformulated in a similar way.

{\it A framework} in $\mbox{\bf R}^n$ is a connected graph
whose vertices are points in  $\mbox{\bf R}^n$ and edges are
straight line segments joining some of its vertices.
It is conventional to call the vertices of a framework
{\it joints} and straight line segments {\it bars}.
The set of 0- and 1-dimensional faces of a polyhedron in $\mbox{\bf R}^n$
can be considered as a typical example of a framework.

A framework in $\mbox{\bf R}^n$ is said to be  {\it flexible}
if it admits a nontrivial analytic deformation, i.e., the positions
of its joints can be changed in $\mbox{\bf R}^n$ analytically with
respect to a parameter in such a way that the length of each bar
remains constant while the distance between some two joints (which
are not joint together by a bar) is not constant.

Sometimes, it is also useful to study so-called {\it pinched frameworks},
i.~e. such frameworks that the spatial positions of some of their joints are fixed and
should not be changed during deformations.

As we have just mentioned above, a suitable framework can be associated
with any polyhedron which consists of the set of its 0- and 1-dimensional faces.
This framework is called a {\it 1-skeleton} of the polyhedron.
Obviously, a polyhedron is flexible if and only if its 1-skeleton
is flexible (the reason is that, according to our definition,
each face of a polyhedron is a simplex).
Hence, the decision problem of whether a given polyhedron is flexible or not
is a particular case of the decision problem of whether a given framework is
flexible or not.
So, we will focus our attention on the latter problem.

The spatial position of a framework is determined as soon as
the positions of its joints $x_i(0)\in \mbox{\bf R}^n$ are given.
The decision problem of whether a given framework is
flexible or not is, obviously, equivalent to the problem of whether
the set of vectors $x_i(0)$ is an isolated solution to (23)
or defines an implicit function in a neighborhood of $x_i(0)$.
Furthermore,  is necessary to exclude trivial deformations
when the framework moves as a rigid body.
This can be done in the same way as it was previously done for polyhedra.
This problem does not appear for pinched frameworks at all, since
they often do not permit any trivial deformations.

In Section 2, we introduced the notion of an approximate solution of degree
$q$ to a system of polynomial equations.
As applied to system (23) associated with a framework,
the term ``infinitesimal flexion of order $q$'' is conventionally used.
More precisely, let a framework in $\mbox{\bf R}^n$ be determined by the
positions $x_i(0)$, $i=1,2,\dots, N$, of its joints.
The set of vectors $x_{i,p}$, $i=1,2,\dots ,N$, $p=1,2,\dots, q$
is said to be an {\it infinitesimal flexion of order $q$} of the framework if
$$\bigg|\sum_{p=1}^{q} x_{i,p} t^p -\sum_{p=1}^{q}
x_{j,p}t^p\bigg|=0\ (\mbox{mod\, } t^q)$$
for all indices $i$ and $j$ such that the joints $x_i(0)$ and $x_j(0)$
are joint by a bar.
An infinitesimal flexion is said to be {\it trivial} if it is an
initial part of the Taylor expansion of the trajectories of $x_i(0)$
under the action of some one-parameter group of isometries of $\mbox{\bf R}^n$.

A framework is said to be {\it $q$th-order infinitesimally flexible} if it
admits a nontrivial infinitesimal flexion of order $q$.
Otherwise it is called {\it $q$th-order infinitesimally rigid}.

Roughly speaking, the following theorem asserts that if a framework admits
a ``regular'' infinitesimal flexion of sufficiently large order then
it is flexible.

{\bf Theorem 7.} {\it Let a framework $P$ be $q$th-order infinitesimally
flexible and let $\sum_{p=0}^{q} X_pt^p$ be a nontrivial infinitesimal
flexion of order $q$.
Let operators $B$ and $C$ be constructed for system (23) which is
associated with $P$ and let there exist a number $k$ \ $(0\leq k < q)$
such that, for every $i=1,2,\dots , q$ and every $j=k,k+1,\dots , q$,
the equation  $$CX=-B(X_i,X_j)-B(X_j,X_i)$$ has a solution lying in the
linear span of the vectors $X_k, X_{k+1}, \dots , X_q$.
Then $P$ is flexible.}

{\bf Proof} of Theorem 7 follows immediately from Theorem 1.

All flexible octahedra in  $\mbox{\bf R}^3$ were classified by R.~Bricard
\cite{Bricard} (see also \cite{Sabitov-E}).
In \cite{Alexandrov}, it is shown (in slightly different
terms) that the conditions of Theorem 7 hold true  with $q=5$ and $k=1$
for the 1-skeleton of so-called Bricard's flexible octahedra of the first type.

The idea to use infinitesimal rigidity (of some order) of a framework
for proving its rigidity is explored for a long time and is based,
first of all, on the following

{\bf Theorem 8.} {\it Every 1-order infinitesimally rigid
framework in  $\mbox{\bf R}^n$ is rigid.}

{\bf Proof} follows immediately from Theorem 2
and from the above reduction of the decision problem of
whether a framework is flexible to the problem of whether a given
solution $x_j=x_j(0)$, $j=1,\dots, N$, to system (23) of
algebraic equations is isolated or the system defines an implicit
function in a neighborhood of the solution.

Theorem 8 is one of the corner stones in problems of existence and
uniqueness for convex polyhedra in the way of exposition which is
used in the classical book \cite{AD}.
Among recent papers that use Theorem 8, we indicate the articles
\cite{Maehara-95} and \cite{Maehara-96} where it is shown that both
in $\mbox{\bf R}^2$ and $\mbox{\bf R}^3$ there exist rigid triangle-free
frameworks with all bars having length 1.
Other statements about interrelations between rigidity and
infinitesimal rigidity may be found in
\cite{Bolker-Roth}, \cite{Con-Serv}, \cite{Maehara-98},
\cite{Whiteley-84TAM}, \cite{Whiteley-84PJM}, \cite{Whiteley-85},
\cite{Whiteley-87}.

The following theorem was proven for the first time
(by different methods) in \cite{Con-Adv}. A similar theorem for
smooth surfaces was obtained in \cite{Efimov-48}.

{\bf Theorem 9.}
{\it Every 2-order infinitesimally rigid
framework in  $\mbox{\bf R}^n$ is rigid.}

{\bf Proof} of Theorem 9 ensues immediately from Theorem 3.

The following theorem generalizes Theorems 8 and 9.

{\bf Theorem 10.} {\it Let a framework $K$ in  $\mbox{\bf R}^n$ have
a single nontrivial linearly independent first-order infinitesimal
flexion and let there exist a number $q\geq 1$ such that $K$ is
$q$th-order infinitesimally rigid. Then $K$ is rigid.}

{\bf Proof.} Eliminate trivial motions as was described above.
The kernel of the operator $C$ which is associated with (23) has
dimension 1.
Put $T= (\mbox{ker\,} C)^\perp$.
According to Theorem 6, (23) has a $T$-standard formal solution
$Y(t)=\sum_{p=0}^{+\infty} Y_pt^p$ such that $Y_0=K$ and $K+Y_1t$
is a nontrivial first-order infinitesimal flexion.
Furthermore, $Y(t)=\sum_{p=0}^{q} Y_pt^p$ is
a nontrivial $q$th-order infinitesimal flexion.
However, this contradicts the hypothesis of Theorem 10
according to which $K$ is $q$th-order infinitesimally rigid.
This contradiction proves Theorem 10.

For smooth surfaces, a theorem similar to Theorem 10 was obtained
(by different methods) in \cite{Sabitov-UGS}.
For some other results about interrelations between higher-order
infinitesimal rigidity and rigidity see, for example,
\cite{Con-Whiteley}, \cite{Stachel-1}, and \cite{Stachel-2}
(for frameworks) and \cite{Perlova-UGS} and \cite{Perlova-MAG}
(for smooth surfaces).

\end{document}